\documentclass[10pt,twoside]{article}
\usepackage{graphicx}
\usepackage{amsmath}
\usepackage{Latex-document}

\title{\bf  Mathematical Results Inspired by Physics\vskip 6mm}

\author{Kefeng Liu\vspace*{-0.5cm}\thanks{Department of Mathematics,
University of California, Los Angeles, CA 90095, USA. E-mail: liu@math.ucla.edu}}
\date{\vspace{-8mm}}

\begin{document}

\maketitle

\thispagestyle{first} \setcounter{page}{457}

\begin{abstract}

\vskip 3mm

I will discuss results of three different types in geometry and
topology. (1) General vanishing and rigidity theorems of elliptic
genera proved by using modular forms, Kac-Moody algebras and
vertex operator algebras. (2) The computations of intersection
numbers of the moduli spaces of flat connections on a Riemann
surface by using heat kernels. (3) The mirror principle about
counting curves in Calabi-Yau and general projective manifolds by
using hypergeometric series.

\vskip 4.5mm


\end{abstract}

\vskip 12mm

\section{Introduction} \label{section 0}\setzero
\vskip-5mm \hspace{5mm }

The results I will discuss are all motivated by the conjectures of
physicists, without which it is hard to imagine that these results
would have appeared. In all these cases the new methods discovered
during the process to prove those conjectures often give us many
more surprising new results. The common feature of the proofs is
that they all depend on localization techniques built upon various
parts of mathematics: modular forms, heat kernels, symplectic
geometry, and various moduli spaces.

Elliptic genera were invented through the joint efforts of
physicists and mathematicians \cite{La}. Actually in Section 2 I
will only discuss in detail a vanishing theorem of the Witten
genus, which is the index of the Dirac operator on loop space.
This is a loop space analogue of the famous Atiyah-Hirzebruch
vanishing theorem. It was discovered in the process of
understanding the Witten rigidity conjectures for elliptic genera.
A loop space analogue of a famous theorem of Lawson-Yau for
non-abelian Lie group actions will also be discussed.

Moduli spaces of flat connections on Riemann surfaces have been
studied for many years in various subjects of mathematics
\cite{AB1}. The computations of the intersection
  numbers on such moduli spaces
have been among the central problems in the subject. In Section 3
I will discuss a very effective way to compute the most
interesting intersection numbers by using the localization
property of heat kernels. This proves several beautiful formulas
conjectured by Witten \cite{W2}. We remark that these intersection
numbers include those needed for the Verlinde formula.

In Section 4 I discuss some remarkable formulas about counting
curves in projective manifolds, in particular in Calabi-Yau
manifolds. I will discuss the mirror principle, a general method
developed in \cite{LLY1}-\cite{LLY4} to compute characteristic
classes and characteristic numbers on
 moduli spaces of stable maps in terms of hypergeometric series.
 The mirror formulas from mirror
  symmetry correspond to the computations of the
 Euler numbers. Mirror principle computes quite general Hirzebruch multiplicative
 classes such as the total Chern classes.

\section{Elliptic genera} \label{section 1}\setzero
\vskip-5mm \hspace{5mm }

 Let $M$ be a
compact smooth spin manifold with a non-trivial $S^1$-action, $D$
be the Dirac operator on $M$. Atiyah and Hirzebruch proved that in
such a situation the index of the Dirac operator ${\mathrm{Ind}}\,
D= \hat{A}(M)=0$, where $\hat{A}(M)$ is the Hirzebruch
$\hat{A}$-genus \cite{AH}. One interesting application of this
result is that a $K3$ surface does not allow any non-trivial
smooth $S^1$-action, because it has non-vanishing $\hat{A}$-genus.

Let $LM$ be the loop space of $M$. $LM$ consists of smooth maps
from $S^1$ to $M$. There is a natural $S^1$-action on $LM$ induced
by the rotation of the loops, whose fixed points are the constant
loops which is $M$ itself. Witten formally applied the
Atiyah-Bott-Segal-Singer fixed point formula to the Dirac operator
on $LM$, from which he derived the following formal elliptic
operator \cite{W}:

$$D^L =D\otimes \bigotimes_{n=1}^\infty S_{q^n}TM= \sum^\infty_{n=0} D\otimes V_n \ q^n$$
 where $q$ is a formal variable and for a vector
 bundle $E$,
 $$ S_q\, E =1+q\, E+\, q^2 \, S^2E +\cdots $$ is the symmetric
 operation and $V_n$ is the combinations of the symmetric products $S^j(TM)$'s by formal
  power series expansion.
So $D^L$, which is called the Dirac operator on loop space,
actually consists of an infinite series of twisted Dirac operators
with the pure Dirac operator $D$ as the degree $0$ term. The index
of $D^L$, denoted by ${\mathrm{Ind}}\, D^L$, is called the Witten
genus. The loop space analogue of our vanishing theorem is the
following:

{\bf Theorem 2.1:} (\cite{Liu}) {\em Let $M$ be a spin manifold
with non-trivial $S^1$-action. Assume $p_1(M)_{S^1} =n\, \pi^*
u^2$ for some integer $n$, then the Witten genus vanishes:
${\mathrm{Ind}}\, D^L =0$.}

Here $p_1(M)_{S^1}$ is the equivariant first Pontrjagin class and
$u$ is the generator of the cohomology group of the classifying
space $BS^1$, and $\pi:\ M\times_{S^1} ES^1 \rightarrow BS^1 $ is
the natural projection from the Borel construction.

 This theorem implies that under the
extra condition on the first Pontrjagin class, we have infinite
number of elliptic operators with vanishing indices. The condition
on the first equivariant Pontrjagin class is equivalent to that
the $S^1$-action preserves the spin structure of $LM$. If we have
a non-abelian Lie group acts on $M$ non-trivially, then for an
$S^1$ subgroup, the condition $p_1(M)_{S^1} =n\, \pi^* u^2$ is
equivalent to $p_1(M)=0$ which implies that $LM$ is spin. As an
easy consequence, we get:

{\bf Corollary 2.2:} {\em Assume a non-abelian Lie group acts on
the spin manifold $M$ non-trivially and $p_1(M)=0$, then the
Witten genus, ${\mathrm{Ind}}\,D^L$, vanishes.}

This corollary should be considered as a loop space analogue of a
result of Lawson-Yau in \cite{LY}, which states that if a
non-abelian Lie group acts on
 the spin manifold $M$ non-trivially, then ${\mathrm{Ind}}\, D=0$. Our results
motivated Hoehn and Stolz to conjecture that, for a compact spin
manifold $M$ with positive Ricci curvature and $p_1(M)=0$, the
Witten genus vanishes. So far all of the known examples have
non-abelian Lie group action, therefore our results applies. It
should be interesting to see how to combine curvature with modular
forms to get vanishing results.

The proof of Theorem 2.1 is an interesting combination of the
Atiyah-Bott-Segal-Singer fixed point formula with Jacobi forms.
The magic combination of geometry and modular invariance implies
the vanishing of the equivariant index of $D^L$. Similar idea can
be used to prove many more rigidity, vanishing and divisibility
results for $D^L$ twisted by bundles constructed from loop group
representations. Such operators can be viewed as twisted Dirac
operators on loop space. See \cite{Liua} and \cite{Liu}. In these
cases the Kac-Weyl character formulas came into play. If we take
the level $1$ representations of the loop group of the spin group
in our general rigidity theorem, we get the Witten conjectures on
the rigidity of elliptic genera \cite{W}, which were proved by
Taubes \cite{T}, Bott-Taubes \cite{BT}, Hirzebruch \cite{H},
Krichever, Landweber-Stong, Ochanine for various cases.

Our method can actually go very far. Recently in \cite{LM} we
proved rigidity and vanishing theorems for families of elliptic
genera and the Witten genus. In \cite{LMZ} we proved similar
theorems for foliated manifolds. In \cite{DLM} such theorems were
generalized to orbifolds. More recently in \cite{DLM1} we have
proved a far general rigidity theorem for $D^L$ twisted by vertex
operator algebra bundles.

If we apply the modular invariance argument to the non-equivariant
elliptic genera, we get a general formula which expresses the
Hirzebruch $L$-form in terms of the twisted $\hat{A}$-forms
\cite{Liu0}. A $12$ dimensional version of this formula, due to
Alveraz-Gaume and Witten, called the miraculous cancellation
formula, had played important role in the development of string
theory. This formula has many interesting mathematical
consequences involving the eta-invariants. We refer the reader to
\cite{Liu0} and \cite{LiuZ}.

\section{Moduli spaces} \label{section 2} \setzero\vskip-5mm \hspace{5mm }

Let $G$ be a compact semi-simple Lie group and ${\cal M}_u$ be the
moduli space of flat connections on a principal flat $G$-bundle
$P$ on a Riemann surface $S$ with boundary, where $u\in Z(G)$ is
an element in the center. Here for simplicity we first discuss the
case when $S$ has one boundary component, $G$ is simply connected
and the moduli space is smooth. A point in ${\cal M}_u$ is an
equivalence class of flat connection on $P$ with holonomy $u$
around the boundary. In general we let ${\cal M}_c$ denote the
moduli space of flat connections on $P$ with holonomy around the
boundary to be $ c\in G$ which is close to $u$, or equivalently in
the conjugacy class of $c$. The following formula is essentially a
refined version of the formula \cite{W2} which Witten derived from
the path integrals on the space of connections.

 {\bf Theorem 3.1:} (\cite{Liu1}, \cite{Liu2}),  {\em We have the following identity:}
$$\int_{{\cal M}_u} p(\sqrt{-1}\Omega)e^{\omega_u}=
|Z(G)|\frac{|G|^{2g-2}}{(2\pi)^{2N_u}}\cdot
\mbox{lim}_{c\rightarrow u}\mbox{lim}_{t\rightarrow
0}\sum_{\lambda\in P_+}
\frac{\chi_\lambda(c)}{d_\lambda^{2g-1}}p(\lambda+\rho)
e^{-tp_c(\lambda)}.$$

The notations in the above formula are as follows: $\omega_u$ is
the canonical symplectic form on ${\cal M}_u$ induced by Poincare
duality on $S$; $p(\sqrt{-1}\Omega)$ is a Pontrjagin class of the
tangent bundle $T{\cal M}_u$ of the moduli space associated to the
symmetric polynomial $p$; $P_+$ is the set of irreducible
representations of $G$ identified as a lattice in $\mathcal{T}^*$
which is the dual Lie algebra of the maximal torus $T$ of $G$;
$p_c(\lambda)=|\lambda+\rho|^2-|\rho|^2$ where
$\rho=\frac{1}{2}\sum_{\alpha\in \Delta^+}\alpha $ with respect to
the Killing form, and $\Delta^+$ denotes the set of positive
roots; $\chi_\lambda$ and $d_\lambda$ are respectively the
character and dimension of $\lambda$; $|G|$ denotes the volume of
$G$ with respect to the bi-invariant metric induced from the
Killing form; $|Z(G)|$ denotes the number of elements in the
center $Z(G)$ of $G$ and finally $N_u$ is the complex dimension of
${\cal M}_u$.

 The starting point for the proof of this theorem is to use the
 holonomy model of the moduli
space and the explicit heat kernel on $G$. We consider the
holonomy map $f: G^{2g} \times O_{c}\rightarrow G$ with $f(x_1,
\cdots, y_g; z)=\prod_{j=1}^g [x_j, y_j]z$ where $O_c$ is the
conjugacy class through the generic point $c\in G$. It is
well-known that the moduli space is given by ${\cal M}_c=
f^{-1}(e)/G$ where $G$ acts on $G^{2g}\times O_{c}$ by
conjugation.

We have the explicit expression for the heat kernel on $G$:
$$H(t, x, y)=\frac{1}{|G|}\sum_{\lambda \in P_+}d_\lambda \cdot
 \chi_\lambda(xy^{-1})e^{-tp_c(\lambda)},$$
where $x,\, y\in G$ are two points. The key idea is to consider
the integral
$$I(t)=\int_{h\in G^{2g}\times O_c} H(t, c, f(h))dh,$$
where $dh$ denotes the induced bi-invariant volume form on $
G^{2g}\times O_c$. We compute $I(t)$ in two different ways. First
as $t\rightarrow 0$, $I(t)$ localizes to an integral on ${\cal
M}_c$, which is the symplectic volume of ${\cal M}_c$ with respect
to the canonical symplectic form induced by the Poincare duality
on the cohomology groups of $S$ with values in the adjoint Lie
algebra bundle. To prove this we used the beautiful observation of
Witten \cite{W1} that the symplectic volume form of ${\cal M}_c$
 is the same as the Reidemeister torsion which arises from the
 Gaussian integral in the heat kernel.

On the other hand the orthogonal relations among the characters of
the representations of $G$ easily give us the infinite sum.
 In summary we have obtained the following more precise version of Witten's
  beautiful formula for the symplectic volume of the moduli space,

{\bf Proposition 3.2:} (\cite{Liu1}) {\em As $t\rightarrow 0$, we
have}

$$\int_{{\cal
M}_c}e^{\omega_c}=|Z(G)|\frac{|G|^{2g-1}|j(c)|}{(2\pi)^{2N_c}|Z_{c}|}\sum_{\lambda\in
P_+}\frac{
\chi_\lambda(c)}{d_\lambda^{2g-1}}e^{-tp_c(\lambda)}+O(e^{-\delta^2/4t}).$$

Here $\delta$ is any small positive number, $|Z_c|$ is the volume
of the centralizer $Z_c$ of $c$, $j(c)=\prod_{\alpha\in
\Delta^+}(e^{\sqrt{-1}\alpha(C)/2}-e^{\sqrt{-1}\alpha(C)/2})$ is
the Weyl denominator, and $N_c$ is the complex dimension of ${\cal
M}_c$.

To get the intersection numbers from the volume formula, we  take
derivatives with respect to $C$ where $c=u\, \mathrm{exp}\,C$.
This is another key observation. By using the relation between the
symplectic form on ${\cal M}_c$ and that on ${\cal M}_u$, and then
taking the limits we arrive at the formula in Theorem 3.1. For the
details see \cite{Liu1} and \cite{Liu2}. Another easy consequence
of the method is that the symplectic volume of ${\cal M}_c$ is a
piecewise polynomial of degree at most $2g\, |\Delta^+|$ in $C\in
{\mathcal T}$ from which we get certain very general vanishing
theorems for those integrals when the degree of the polynomial $p$
is relatively large \cite{Liu2}.

Similar results for moduli spaces when $S$ has more boundary
components can be obtained in the same way \cite{Liu2}. More
precisely, assume $S$ has $s$ boundary components and consider the
moduli space of flat connections on the principal $G$ bundle $P$
with holonomy $c_1, \cdots ,c_s\in G$ around the corresponding
boundaries. Let ${\cal M}_{\bf c}$ denote the moduli space and
$\omega_{\bf c}$ denote the canonical symplectic form. Then we
have

{\bf Theorem 3.3:} (\cite{Liu2}) {\em The following formula holds:}
$$ \int_{{\cal{M}}_{{\bf c}}} p(\sqrt{-1}\Omega)e^{\omega_{{\bf c}}} =
|Z(G)|\frac{|G|^{2g-2+s}\prod_{j=1}^s \limits j(c_j)}{(2\pi)^{2N_{\bf c}}{\prod_{j=1}^s \limits |Z_{c_j}|}} \
\lim_{t\rightarrow 0} \sum_{\lambda\in P_+} \frac{\prod_{j=1}^s \limits
\chi_\lambda(c_j)}{d_\lambda^{2g-2+s}}p(\lambda+\rho)e^{-tp_c(\lambda)}. $$

Here $N_{\bf c}$ is the complex dimension of ${\cal M}_{\bf c}$
and $p(\sqrt{-1}\Omega)$ is a Pontrjagin class of ${\cal{M}}_{{\bf
c}}$. By taking derivatives with respect to the $c_j$'s we can get
intersection numbers involving the other generators of the
cohomology ring of ${\cal M}_{\bf c}$, as well as the polynomial
property. From index formula we know that the integrals in our
formulas contain all the information needed for the famous
Verlinde formula. Recently the general Verlinde formula has been
directly derived along this line of idea \cite{BL}.

This localization method of using heat kernels can be applied to
other general situation like moment maps, from which we derive the
non-abelian localization formula of Witten. See \cite{Liu2} for
applications to three dimensional manifolds and see \cite{Liu3}
for applications involving finite groups and moment maps.

\section{Mirror principle} \label{section 3} \setzero\vskip-5mm \hspace{5mm }

 Let $X$ be a projective manifold. Let ${\cal M}_{g, k}(d, X)$ denote the
moduli space of stable maps of genus $g$ and degree $d$ with $k$
marked points into $X$. Modulo the obvious equivalence, a point in
${\cal M}_{g, k}(d, X)$ is given by a triple $ (f; C;\, x_1,
\cdots, x_k)$ where $ f: C\rightarrow X$ is a degree $d$
holomorphic map and $x_1, \cdots, x_k$ are $k$ points on the genus
$g$ curve $C$. Here $d\in H_2(X, ,{\mathbf Z})$ will also be
identified as the integral index $(d_1, \cdots, d_n)$ by choosing
a basis of $H_2(X, ,{\mathbf Z})$ dual to a basis of Kahler
classes.

This moduli space may have higher dimension than expected. Even
worse, its different components may have different dimensions. To
define integrals on such space, we need the virtual fundamental
cycle first constructed in \cite{LT} and later in \cite{BF}. Let
us denote by $LT_{g, k}(d,X)$ the virtual fundamental
 cycle which is a homology class of the expected dimension in ${\cal M}_{g, k}(d,
 X)$.

We first consider the case $k=0$. Let $V$ be a concavex bundle on
$X$. The notion of concavex bundles was introduced in \cite{LLY1},
it is a direct sum of a positive and a negative bundle on $X$.
From a concavex bundle $V$, we can obtain a sequence of vector
bundles $V^g_d$ on $M_{g, k}(d, X)$ by taking either $H^0(C,
f^*V)$ or $H^1(C, f^*V)$, or their direct sum. Let $b$ be a
multiplicative characteristic class. The main problem of mirror
principle is to compute the integral \cite{Ko}

$$K^g_d= \int_{LT_{g, 0}(d, X)} b(V^g_d).$$
More precisely, let $\lambda$ and $T=(T_1,\cdots, ,T_n)$ be formal
variables. Mirror principle is to compute the generating series,

$$F(q, \lambda)=\sum_{d,\, g} K^g_d \,\lambda^g ~e^{d\cdot T}$$ in terms of certain natural explicit
hypergeometric series. So far we have rather complete picture for
the case of balloon manifolds.

A balloon manifold $X$ is a projective manifold with complex torus
action and isolated fixed points. Let $H=(H_1, \cdots, H_n)$ be a
basis of equivariant Kahler classes. Then $X$ is called a balloon
manifold if $H(p)\neq H(q)$ when restricted to any two fixed
points $p,\, q\in X$, and the tangent bundle $T_pX$ has linearly
independent weights for any fixed point $p\in X$. The complex
$1$-dimensional orbits in $X$ joining every two fixed points in
$X$ are called balloons which are copies of ${\mathbf P}^1$. We
require the bundle $V$ to have fixed splitting type when
restricted to each balloon \cite{LLY1}.

{\bf Theorem 4.1:} (\cite{LLY1}-\cite{LLY4}) {\em  Mirror
principle holds for balloon manifolds and concavex bundles.}

In the most interesting cases for the mirror formulas, we simply
take characteristic class $b$ to be the Euler class and the genus
$g=0$. The mirror principle implies that mirror formulas actually
hold for very general manifolds such as Calabi-Yau complete
intersections in toric manifolds and in compact homogeneous
manifolds. See \cite{LLLY} and \cite{LLY3} for details. In
particular this implies all of the mirror formulas for counting
rational curves predicted by string theorists. Actually mirror
principle holds even for non-Calabi-Yau and for certain local
complete intersections. In \cite{LLY4} we developed the mirror
principle for counting higher genus curves, for which the only
remaining problem is to find the explicit hypergeometric series.
Also our method clearly works well for orbifolds.

As an example, we consider a toric manifold $X$ and genus $g=0$.
Let $D_1,..,D_N$ be the toric invariant divisors, and $V$ be the
direct sum of line bundles: $V=\bigoplus_j L_j$ with
$c_1(L_j)\geq0$ and $c_1(X)=c_1(V)$. We denote by $\langle\cdot,
\cdot\rangle$ the pairing of homology and cohomology classes. Let
$b$ be the Euler class and
$$\Phi(T)=\sum_d K^0_d \,e^{d\cdot T}$$ where $d\cdot T=d_1T_1+\cdots d_nT_n$. Introduce the
hypergeometric series
$$
HG[B](t)=e^{-H\cdot t}\sum_d \prod_j\prod_{k=0}^{\langle
c_1(L_j),d\rangle} (c_1(L_j)-k){\prod_{\langle D_a,d\rangle<0}
\prod_{k=0}^{-\langle D_a,d\rangle-1}(D_a+k)\over \prod_{\langle
D_a,d\rangle\geq0} \prod_{k=1}^{\langle D_a,d\rangle}(D_a-k)}~
e^{d\cdot t}
$$ with $t=(t_1,\cdots, t_n)$ formal variable.

{\bf Corollary 4.2:} \cite{LLY3} {\em There are explicitly
computable functions $f(t),\, g(t)=(g_1(t), \cdots, g_n(t))$, such
that
$$
\int_X\left(e^f HG[B](t)-e^{-H\cdot T}e(V)\right) =2\Phi-\sum_j
T_j {\partial\Phi\over\partial T_j}
$$
where $T=t+g(t)$.}

From this formula we can determine $\Phi(T)$ uniquely. The
functions $f$ and $g$ are given by the expansion of $HG[B](t)$. We
can also replace $V$ by general concavex bundles \cite{LLY3}. To
make our algorithm more explicit, let us consider the the
Calabi-Yau quintic, for which we have the famous Candelas formula
\cite{Ca}. In this case $V=\mathcal{O}(5)$ on $X={\mathbf P}^4$
and the hypergeometric series is:

$$HG[B](t)= e^{H\,t}\sum^\infty_{d=0}
\frac{\prod_{m=0}^{5d}(5H+m)}{\prod_{m=1}^{d}(H+m)^5}\, e^{d\,t},
$$where $H$ is the hyperplane class on $\mathbf{P}^4$ and $t$ is a parameter.
Introduce

$${\cal{F}}(T) =\frac{5}{6}T^3 +\sum_{d>0} K^0_d\, e^{d\,T}.$$
The algorithm is to take the expansion in $H$:

$$HG[B](t)= H\{ f_0(t) +f_1(t)H +f_2(t) H^2 +f_3(t)H^3\}.$$
 Then the famous Candelas formula can be reformulated as

 {\bf Corollary 4.3:} (\cite{LLY1}) {\em With $T=f_1/f_0$, we have}
 $${\cal F}(T)=
\frac{5}{2}(\frac{f_1}{f_0}\frac{f_2}{f_0}-\frac{f_3}{f_0}).$$

 Another rather interesting consequence of mirror principle is the
local mirror symmetry which is the case when $V$ is a concave
bundle. Local mirror symmetry is called geometric engineering in
string theory which is used to explain the stringy origin of the
Seiberg-Witten theory. In these cases the hypergeometric series
are the periods of elliptic curves which are called the
Seiberg-Witten curves. These elliptic curves are the mirror
manifolds of the open Calabi-Yau manifolds appeared in the local
mirror formulas. For example the total space of the canonical
bundle of a del Pezzo surface is an example of open Calabi-Yau
manifold covered by the local mirror symmetry. The case ${\mathbf
P}^2$ already has drawn a lot of interests in string theory. The
case for $\mathcal{O}(-1)\oplus \mathcal{O}(-1)$ on ${\mathbf
P}^1$ easily gives the multiple cover formula.

 The key ingredients for the proof of the mirror principle consist of the
following: linear and non-linear sigma model, Euler data, balloon
and hypergeometric Euler data. As explained in \cite{LLY4}, these
ingredients are independent of the genus of the curves, except the
hypergeometric Euler data, which for $g>0$ is more difficult to
find out, while for the genus $0$ case it can be easily read out
from localization at the smooth fixed points of the moduli spaces
which are covers of the balloons. The interested reader is refered
to \cite{LLY1}-\cite{LLY4} for details. Our idea is to go to the
equivariant setting and to use the localization formula as given
in \cite{AB} and its virtual version in \cite{GP} on two moduli
spaces which we called non-linear and linear sigma models. One key
observation is the functorial localization formula
\cite{LLY1}-\cite{LLY4}. We apply this formula to the equivariant
collapsing map between the two sigma models, and to the evaluation
maps. One can see \cite{LLY4} for the existence of the collapsing
map for arbitrary genus. Hypergeometric series naturally appear
from localizations on the linear sigma models and at the smooth
fixed points in the moduli spaces.

Euler data is a very general notion, it can include general
Gromov-Witten invariants by adding the pull-back classes by the
evaluation map $ev_j$ at the marked points. More precisely we can
try to compute integrals of the form:
$$K^g_{d, k}=\int_{LT_{g, k}(d, X)}\prod_{j}ev_j^*\omega_j
\,\cdot  b(V^g_d)$$ where $\omega_j\in H^*(X)$. By introducing the
generating series with summation over $k$, we can still get Euler
data. One goal of the most general mirror principle is to
explicitly compute such series in terms of hypergeometric series.

We remark that the development of the proof of the mirror formulas
owes to many people, first to the string theorists Candelas and
his collaborators, Witten, Vafa, Warner, Greene, Morrison, Plesser
and many others. They used the physical theory of mirror symmetry,
and their computations used mirror manifolds and their periods.
For the general theory of mirror principle, see
\cite{LLY1}-\cite{LLY4}. See also \cite{Gi}, \cite{Ga} and
\cite{B} for different approaches to the mirror formula.

\section{Concluding remarks} \label{section 5} \setzero\vskip-5mm \hspace{5mm }

 Localization techniques have been very
successful in solving many conjectures from physics. In the
meantime string theorists have produced many more exciting
 new conjectures. We can certainly expect their solutions by using
localizations.

Recently several mirror formulas of counting holomorphic discs
have been conjectured by Vafa and his collaborators. The boundary
of the disc is mapped into a Langrangian submanifolds of the
Calabi-Yau. Other related conjectures include the Gopakumar-Vafa
conjecture on the higher genus multiple covering formula and the
mirror formulas for counting higher genus curves in Calabi-Yau.
With Chien-Hao Liu we are trying to extend the mirror principle to
these settings. Another exciting conjecture is the S-duality
conjecture which includes the Witten conjecture on the equality of
the Donaldson invariants with the Seiberg-Witten invariants and
the Vafa-Witten conjecture on the modularity of the generating
series of the Euler numbers of the moduli spaces of self-dual
connections. Some progresses are made by constructing a larger
moduli space with circle action, the so-called non-abelian monople
moduli spaces. Finally there is the
Dijkgraaf-Moore-Verlinde-Verlinde conjecture on the generating
series of the elliptic genera of Hilbert schemes. For an approach
of using localization, see \cite{LMZh}.

\label{lastpage}

\end{document}